\documentclass[a4paper,12pt]{article}
\usepackage[left=2cm,right=2cm,top=2cm,bottom=2cm,bindingoffset=0cm]{geometry}
\usepackage[T1]{fontenc}
\usepackage[utf8]{inputenc}
\usepackage{latexsym}
\usepackage{amsmath}
\usepackage{amssymb}
\usepackage[english,main=ukrainian]{babel}
\usepackage{bm}
\usepackage{graphicx}
\usepackage{wrapfig}
\usepackage{fancybox}
\usepackage{amsthm}
\usepackage{amsthm}

\newcounter{punkt}[section]
\newcounter{theorem}[section]
\newcounter{lemma}[section]
\newcounter{definite}[section]
\newcounter{example}[section]
\newcounter{tv}[section]
\newcommand{\theor}{\par\refstepcounter{theorem}%
{\bf Теорема \thetheorem .}\,\,}

\newcommand{\defi}{\par\refstepcounter{definite}%
{\bf Означення \thedefinite .}\,\,}
\newcommand{\exam}{\par\refstepcounter{example}%
{\bf Приклад \theexample .}\,\,}
\newcommand{\tv}{\par\refstepcounter{tv}%
{\bf Твердження \thetv .}\,\,}

\author{
  Максим Ткачук\\
  \texttt{maxim.v.tkachuk@gmail.com}
}

\title { Пошарова опуклість. }

\begin{document}

\maketitle

Класичним визначенням опуклості в афінному просторі над полем
дійсних чисел є наступне: множина є опуклою, якщо разом із двома
своїми точками вона містить і відрізок, що їх з'єднує. Іншими
словами множина опукла, якщо її перетин з довільною прямою зв'язний.
Поруч із цим визначенням існує поняття лінійної опуклості.
Підмножина афінного простору називається лінійно опуклою, якщо через
будь-яку точку доповнення до неї проходить гіперплощина, яка дану
множину не перетинає. Іншими словами, множина лінійно опукла, якщо
доповнення до неї є об'єднанням гіперплощин. В афінному просторі
опуклі відкриті множини (опуклі компакти) це в точності всі зв'язні
лінійно опуклі відкриті множини (лінійно опуклі компакти), а всяка
зв'язна компонента лінійно опуклої відкритої множини опукла. Також
кожна регулярна і обмежена лінійно опукла область є опуклою. Проте
при відкиданні умови обмеженості або регулярності дане твердження
перестає бути вірним, що показують наступні приклади:

1. регулярної але необмеженої лінійно опуклої області

\begin{figure}[h!]
  \centering
  \includegraphics[width=150pt]{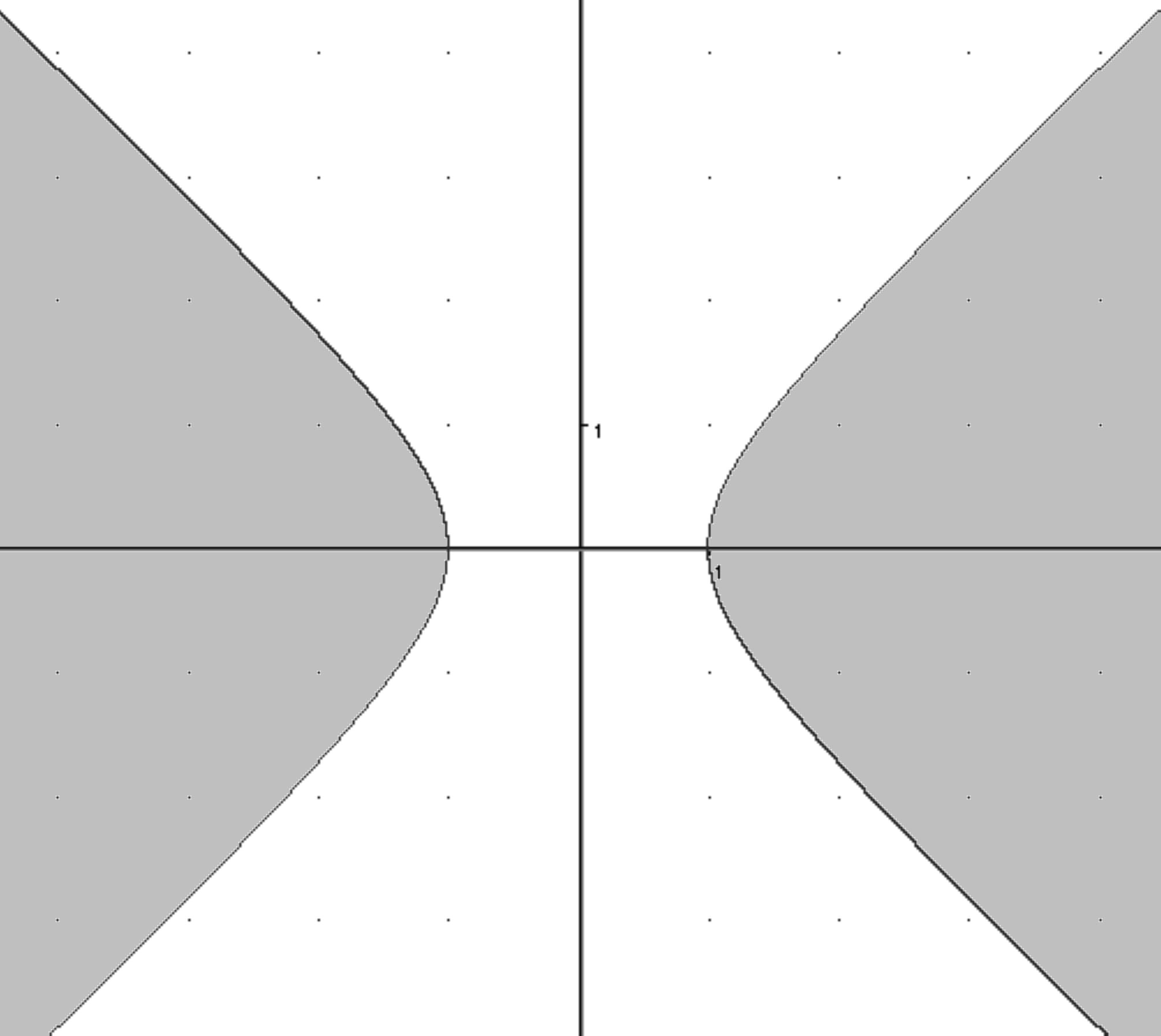}
\end{figure}


2. обмеженої але нерегулярної лінійно опуклої області

\begin{figure}[h!]
  \centering
  \includegraphics[width=300pt]{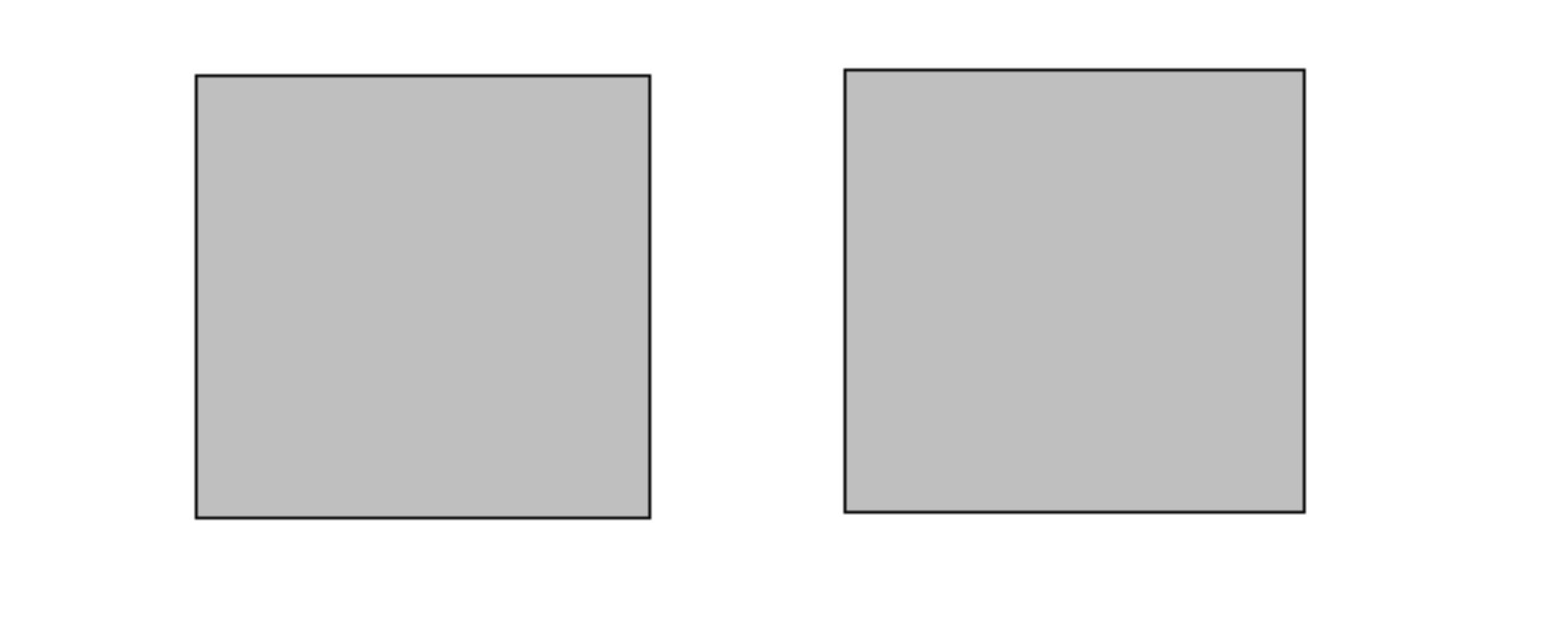}
\end{figure}

 Узагальненням даних
понять є визначення опуклості й лінійної опуклості в комплексному
$n$-вимірному евклідовому просторі $\mathbb{C}^n$. Тобто, множина
називається $\mathbb{C}$-опуклою, якщо її перетин довільною
комплексною прямою ациклічний (зв'язний й однозв'язний). Множина
називається лінійно $\mathbb{C}$-опуклою, якщо її доповнення до
$\mathbb{C}^n$ є об'єднанням комплексних гіперплощин. Також існує
поняття m-опуклості, де замість прямих і гіперплощин беруться
m-вимірні площини. Поняття лінійної $\mathbb{C}$-опуклості в
$\mathbb{C}^n$ аналогічне лінійній $(2n-2)$-опуклості в
$\mathbb{R}^n$, тільки замість усіх $(2n-2)$-площин беруться тільки
ті, які відповідають комплексним гіперплощинам. З огляду на це,
можна ввести поняття лінійної опуклості в $\mathbb{R}^{n}$ і
відносно деякої множини $W$ афінних підпросторів в $\mathbb{R}^{n}$

 \vskip 0.5cm

\defi {\em Множина $E\subset \mathbb{R}^{n}$
називається опуклою відносно множини $W$ афінних лінійних
підмноговидів в $\mathbb{R}^{n}$, якщо її доповнення до
$\mathbb{R}^{n}$ є об'єднанням елементів з $W$.}

\defi {\em Для довільної множини $E\subset
\mathbb{R}^{n}$ назвемо підмножину $E^{*}\subset W$ спряженою до
множини $E$, якщо $\displaystyle E^{*}=\{l\in W \mid l$ не
перетинає $E\},$ відповідно $E^{**}=\mathbb{R}^{n}\backslash
\bigcup \limits_{{l\in E^{*}}} l\subset \mathbb{R}^{n}$ (остання
формула коректна, оскільки кожен елемент множини $W$ є підмножиною
$\mathbb{R}^{n}$).}

Помітимо, що поняття спряженої і подвійної спряженої множини має
сенс тільки для підмножин $\mathbb{R}^{n}$ оскільки множина $W$
може й не мати структури лінійного простору. Очевидно, що
опуклість множини $Е$ відносно $W$ рівносильна рівності
$E=E^{**}$. Множину $E ^{**}$ назвемо опуклою оболонкою множини $E
\in \mathbb{R}^{n}$ відносно $W$. Без усяких додаткових припущень
відносно множини $W$ справедливі наступні твердження.

 \vskip 0.5cm
\tv \label{tv1} {\em Якщо множина $E_{1}$ належить множині $E$, то
для спряжених множин справедливі наступні включення $E^{*}\subset
E^{*}_{1}$ й $E_{1}^{**}\subset E^{**}$.}

\textbf{Доведення.} Якщо пряма $l$ не перетинає $E$, то вона не
перетинає й $E_{1}$, отже $E^{*}$ належить $E^{*}_{1}$. $ \bigcup
\limits_{{l\in E^{*}}} l \subset \bigcup \limits_{{l\in
E_{1}^{*}}} l$ а отже,
$$E_{1}^{**}=\mathbb{R}^{n}\backslash
\bigcup \limits_{{l\in E_{1}^{*}}} l\subset
\mathbb{R}^{n}\backslash \bigcup \limits_{{l\in E^{*}}} l=E^{**}.
\blacksquare$$

\tv {\em Нехай $f: \mathbb{R}^{n} \rightarrow \mathbb{R}^{m} $ -
сюр'єктивне афінне відображення й множина $E \subset
\mathbb{R}^{m}$ опукла відносно множини $W$, тоді множина
$f^{-1}(E) \subset \mathbb{R}^{n}$ опукла відносно множини
$f^{-1}(W)= \{ f^{-1}(p)\subset \mathbb{R}^{n}, p \in W \}$.}

\textbf{Доведення.} Візьмемо довільну точку $x$, яка не належить
прообразу $f^{-1}(E)$, тоді існує афінний підпростір $l \in W$, що
проходить через точку $f(x)$ і не перетинає множину $E$. Отже й
прообрази $f^{-1}(l)$ і $f^{-1}(E)$ не перетинаються. Тепер для
повноти доведення варто лише уточнити, що $f^{-1}(l)$ є афінним
підпростором $\mathbb{R}^{n}$.$\blacksquare$

\tv \label{tv3} {\em  $(\bigcup
\limits_{\alpha}E_{\alpha})^{*}=\bigcap
\limits_{\alpha}E_{\alpha}^{*}$.}

\textbf{Доведення.} Розглянемо ланцюжок рівностей: $\displaystyle
(\bigcup \limits_{\alpha}E_{\alpha})^{*}= \{l\in W \mid l$ не
перетинає $\bigcup \limits_{\alpha}E_{\alpha}\}=\bigcap
\limits_{\alpha}\{l\in W \mid l$ не перетинає $E_{\alpha}\}
 =\bigcap \limits_{\alpha} E_{\alpha}^{*}$. $\blacksquare$

\tv {\em  $E\subset E^{**}$, $(E^{**})^{*} = E^{*}$.}

\textbf{Доведення.} Об'єднання всіх лінійних многовидів з множини
$W$, які не перетинають $E$, належить доповненню до $E$: $\bigcup
\limits_{l\bigcap E= \emptyset, l \in W} l \subset
\mathbb{R}^{n}\backslash E$. Перейшовши до доповнень отримаємо
включення $E\subset E^{**}$. Із твердження \ref{tv1} слідує, що
$E^{***} = (E^{**})^{*} \subset E^{*}$. Далі, нехай лінійний
многовид $l$ є елементом множини $E^{*}$, тоді $l$ не перетинає
множину $E$, а отже $l$ не перетинає й множину
$E^{**}=\mathbb{R}^{n}\backslash \bigcup \limits_{{l\in E^{*}}}
l$, звідки $l \in (E^{**})^{*}$. Отримали обернене включення
$E^{*} \subset (E^{**})^{*}= E^{***}$.$\blacksquare$

\tv {\em  Перетин довільної сукупності опуклих відносно $W$ множин
опуклий відносно $W$ (аксіома опуклості).}

\textbf{Доведення.} Розглянемо перетин $E=\bigcap E_{k}$, де кожна
множина $E_{k}$- опукла відносно $W$. Виберемо деяку точку $x$
поза множиною $E$, тоді $x\notin E_{k}$ для деякого номера $k$.
Тоді існує афінний многовид $l \in W$, який містить точку $x$ і
для якого перетин $l\bigcap E_{k}$ порожній, а значить й $l\bigcap
E = \emptyset$.$\blacksquare$

Для об'єднання множин аналогічний результат невірний навіть для
послідовностей множин впорядкованих по включенню.

\exam  Розглянемо послідовність множин на площині, опуклих
відносно множини всіх прямих (лінійно опуклих)
$$D_n = \{(x,y) | x \geq 0, -nx \leq y \leq nx\}.$$
Їх об'єднання, множина
$$D = \bigcup \limits_n D_n = \{(x,y) | x > 0\}\bigcup\{(0,0)\}$$
вже не є лінійно опуклою. Кожна з прямих, що проходять через точку
$(0,1)$ перетинає множину $D$.

 Проте для компактних множин є справедливими наступні твердження.

\tv {\em  Нехай задана монотонна послідовність компактів $E_{k}$,
$k=1,2,..$, $E_{k+1} \subset E_{k}$  й $E= \bigcap \limits
_{k}E_{k}^{}$. Тоді $E^{*}= \bigcup \limits _{k}E_{k}^{*}$.}

\textbf{Доведення.} Множина $E$ належить $E_{k}$, тоді для всіх
номерів $k$ з твердження \ref{tv1} випливає включення
$E_{k}^{*}\subset E^{*}$, і отже $\bigcup \limits_{k}
E_{k}^{*}\subset E^{*}$ .
 Доведемо включення в іншу сторону. Припустимо, що перетин $l\bigcap E$ є
 порожньою множиною,  але для всіх номерів $k$ перетин
 $l\bigcap E_{k} =P_{k}\neq \emptyset$.
 Всі множини $P_{k}$ - компактні і справедливе включення
 $P_{k+1} = l\bigcap E_{k+1} \subset l\bigcap E_{k} = P_{k}$.
 Але тоді з теореми про вкладені компакти випливає, що перетин
 $l\bigcap E=\bigcap \limits _{k}P_{k}$ не є порожньою множиною.
 Отримали протиріччя, а значить, $E^{*}= \bigcup \limits
_{k}E_{k}^{*}$ .$\blacksquare$

\tv \label{tv7}{\em Якщо компакт $Е$ апроксимується ззовні
послідовністю областей $D_{k}, k=1,2,..$ , то $E^{*}= \bigcup
\limits _{k}D_{k}^{*}$, $E^{**}= \bigcap \limits _{k}D_{k}^{**}$.}

\textbf{Доведення.} Із твердження \ref{tv1} випливає включення
$\bigcup \limits _{k} D_{k}^{*} \subset E^{*}$. Якщо афінний
многовид $l \in W$ не перетинає $E$, то внаслідок компактності
множини $E$, $l$ не перетинає й деякий окіл $E$, а внаслідок
апроксимації $E$ областями $D_{k}, k=1,2,..$, цей окіл містить і
деяку область $D_{k_0}$, отже $E^{*} \subset D_{k_0}^{*} \subset
\bigcup \limits _{k}D_{k}^{*}$. Рівність $E^{**}= \bigcap \limits
_{k}D_{k}^{**}$ випливає із твердження \ref{tv3}. $\blacksquare$

\defi {\em Відкрита множина $D\subset \mathbb{R}^{n}$
називається cлабо опуклою відносно множини $W$ афінних лінійних
підмноговидів в $\mathbb{R}^{n}$, якщо через кожну точку її межі
$\partial D$ проходить деякий елемент з множини $W$.}

Очевидно будь-яка опукла відносно множини $W$ область є cлабо
опуклою відносно  $W$.

\tv \label{tv8} {\em Кожна слабо  опукла відносно  $W$ область $D$
є зв'язною компонентою деякої опуклої відносно  $W$ відкритої
множини. Крім того, кожна зв'язна компонента опуклої відносно $W$
відкритої множини є слабо опуклою відносно  $W$.}

\textbf{Доведення.} Нехай   $x$ - довільна точка межі $\partial
D$. Виберемо один із елементів $l_x$ множини $W$, що проходять
через $x$ та не перетинають $D$.  Отже, кожна точка, у тому числі
й $x$, площини $l_x$  не належить $D^{**}$. Тобто, будь-яка точка
межі $\partial D$ не належить $D^{**}$. Але $D \subset D^{**}$,
тому $D$ - зв'язна компонента опуклої відносно  $W$  множини
$D^{**}$. Другий висновок очевидний.

\defi {\em  Замкнену множину $M\subset \mathbb{R}^{n}$
  назвемо слабо опуклою відносно $W$,
якщо вона апроксимується ззовні послідовністю слабо опуклих
відкритих множин.}

\tv  {\em Кожний слабо   опуклий  зв'язний компакт $K$ є зв'язною
компонентою деякої опуклої замкненої множини.}

\textbf{Доведення.} Нехай $K$ апроксимується ззовні послідовністю
слабо  опуклих областей   $D_k$, $K=\bigcap \limits_k D_k$.
Відповідно до твердження \ref{tv8} $D_k$ співпадає із зв'язною
компонентою множини $D_k^{**}$, а  згідно твердження \ref{tv7}
$K^{**}=(\bigcap \limits_k D_k)^{**}=\bigcap \limits_k D_k^{**}$.
Тоді зв'язна компонента $\bigcap \limits_k D_k^{**}$, яка містить
$K$,  міститься в перетині $\bigcap \limits_k D_k=K$, обернене
включення очевидне.

 \textbf{Зауваження.} Замість простору
$\mathbb{R}^{n}$ можна взяти довільний афінний многовид.

При деяких додаткових припущеннях, а саме при локальній
тривіальності введенного в самому тексті доведення розшарування
$\Sigma$, справедливе наступне твердження.

 \theor  {\em  Область $D \subset \mathbb{C}^n$ з гладкою межею
 класу $C^1$ тоді і тільки тоді є лінійно опуклою, коли існує
 неперервна функція $f$, яка ставить у відповідність кожній
 комплексній гіперплощині $L$, що перетинає область $D$, точку в перетині
 $L\bigcap D$. }

 \textbf{Доведення.}  Спочатку доведемо достатність.

 Припустимо, що існує вказана в умові теореми неперервна функція $f$, але
  $D$ не є лінійно опуклою, тоді лінійна оболонка області  $D$ не
  співпадає з самою областю. Тобто існує така $z$, що не належить
  $D$і належить її лінійній оболонці. Довільна комплексна
  гіперплощина, що проходить через точку $z$, перетинає область
  $D$.

  Кожному одиничному вектору $\nu \in \mathbb{C}^n$ відповідає
  гіперплощина $L_\nu \ni z$, для якої $\nu$ є нормаллю. Внаслідок
  припущення, кожній $L_\nu$ поставлена точка $z_\nu \in L_\nu \bigcap
  D$, $z_\nu \neq z$, оскільки $z \notin D$.

  Поставимо у відповідність кожному вектору $zz_\nu$ колінеарний
  йому (з дійсним коефіцієнтом) одиничний вектор $\tau$.

  Ми побудували неперервну функцію $\varphi$, оскільки вона є
  композицією неперервних функцій, яка ставить довільному
  одиничному вектору $\nu$ одиничний вектор $\tau$

  $$\varphi : S^{2n-1} \rightarrow S^{2n-1}.$$ Але за схемою побудови
  $\varphi (-z) = \varphi (z)$, $\varphi (z) \perp z$, що
  неможливо, оскільки з ортогональності   $\varphi (z) \perp z$
  випливає гомотопність відображення $\varphi$ тотожньому
  відображенню $id : S^{2n-1} \rightarrow S^{2n-1}$. Отже $\varphi$
  має степінь 1 по модулю 2, а з умови $\varphi (-z) = \varphi
  (z)$ випливає, що $\varphi$ має степінь $0 (mod \ 2)$.

 Необхідність.

 Нехай $D$ - лінійно опукла область з гладкою межею. Звідси
 випливає, що вона є також сильно лінійно опуклою. Існує
 гомеоморфізм між множиною гіперплощин, що перетинають $D$ і множиною
 $\mathbb{C}^n \backslash D^* $, доповнення до спряженої множини $D^*$.
  Гіперплощину, що відповідає точці $z \in \mathbb{C}^n \backslash
  D^*$ позначимо через $L_z$. Внаслідок сильної лінійної опуклості
  $D$, множина $\mathbb{C}^n \backslash D^* $  є ациклічною і
  відкритою, і кожній точці  $z \in \mathbb{C}^n \backslash D^* $
  відповідає ациклічна множина $F(z) = L_z \bigcap D$.

  Доведемо, що багатозначне відображення $F$  є неперервним.
  Тобто, кожен окіл множини $F(z_0)$  містить всі $F(z)$ для всіх $z$
  з деякого околу точки $z_o$ (неперервність зверху) і довільний окіл
  довільної точки з $F(z_o)$  має непорожній перетин з усіма $F(z)$
  для всіх $z$ з деякою околу точки $z_0$ (неперервність знизу).

  Довільна гіперплощина $L_z$, що перетинає $D$, не може бути
  дотичною до $\partial D$ в жодній точці перетину, оскільки
  область $D$ є лінійно опуклою, а отже і локально лінійно
  опуклою, а отже і локально лінійно опуклою і
  кожна дотична до $\partial D$
  гіперплощина не повинна перетинати $D$.

  З останнього твердження і випливає неперервність відображення
  $F$, а також, що множини $L_z \bigcap D$ мають гладку межу.

  Оскільки всі множини  $L_z \bigcap D$ ациклічні, відкриті і
  мають гладку межу, то вони попарно гомеоморфні. Тоді ми можемо
  розглянути відображення  як розшарування $\Sigma$ над базою
  $\mathbb{C}^n \backslash D^*$
  і з шаром  $L_z \bigcap D$. Воно є локально тривіальним за припущенням.
  А внаслідок стягуваності бази
  $\mathbb{C}^n \backslash D^*$ воно є топологічно еквівалентним
  декартовому добутку, при умові, що саме розшарування ми
  розглядаємо у просторі $\mathbb{C}^n \backslash D^* \times D$.
  Вибравши довільним чином точку в одному шарі декартового добутку
  $w \in L_z \bigcap D$ ми отримаємо неперервну однозначну функцію
  $f$ таку, що $f(z) \in F(z)$, для всіх $z \in \mathbb{C}^n \backslash
  D^*$, що і треба було довести.

\vskip 0.5cm
\setcounter{theorem}{0} \setcounter{lemma}{0}
\addtocounter{section}{1} \setcounter{punkt}{0}
\section*{\thesection. \
Деякі варіанти різних означень лінійної опуклості}
\addcontentsline{toc}{section}{\thesection. \ Деякі варіанти
різних означень лінійної опуклості}

Наведемо декілька прикладів, коли спряжений простір ізоморфний
даному, і  кожній точці даного простору відповідає лінійний
підмноговид спряженого простору, що складається із всіх елементів,
що проходять через цю точку. У цьому випадку $E^{**}=(E^{*})^{*}$,
де в правій частині рівності обидва спряження беруться по
однакових формулах але в різних просторах, а $E$ належить одному з
них.

\vskip 5mm
\exam  Розглянемо простір $\mathbb{R}^{3}$ й $W$ - множину прямих
у ньому, паралельних даній площині $S$. Тоді опуклість множини $Е$
відносно $W$ рівносильна класичній лінійній опуклості кожного
перетину множини $Е$ площиною, паралельною площині $S$. Якщо ми
введемо в $\mathbb{R}^{3}$ координати так, що
$S=\{x=(x_{1},x_{2},x_{3})\in \mathbb{R}^{3} | x_{3}=0 \}$, тоді
кожна пряма з $W$ буде задаватися рівнянням
$$l=\{x=(x_{1},x_{2},x_{3})\in \mathbb{R}^{3} | x_{3}=d,
ax_{1}+bx_{2}+cx_{3}=0,$$ $$(a:b:c)\in \mathbb{R}P^{2}, d\in
\mathbb{R} \}.$$ А значить множина $W$ гомеоморфна декартовому
добутку $\mathbb{R}P^{2}\times \mathbb{R}$ , де $\mathbb{R}P^{2}$
- проективна площина.
  Якщо ми розширимо $\mathbb{R}^{3}$ до простору
  $\mathbb{R}P^{2}\times \mathbb{R}$ й введемо координати
  $x=(x_{1}:x_{2}:x_{3};x_{4})$, де
  $(x_{1}:x_{2}:x_{3})$ проективні координати в проективній площині
   $\mathbb{R}P^{2}$, $x_{4}\in \mathbb{R}$,
  тоді спряження буде задаватися у вигляді, симетричному відносно даного
  й спряженого просторів.
   Отже, дані відображення індукують на множині $W$ структуру
   лінійного простору $\mathbb{R}P^{2}\times \mathbb{R}$, і ми
   можемо розглядати множини $E$ і $E^{*}$ як такі, що належать
   одному лінійному простору $\mathbb{R}P^{2}\times \mathbb{R}$.

\vskip 5mm

\exam Візьмемо простір $\mathbb{R}^{3}$ і множину $W$ прямих,  що
належить площинам пучка,
 кожна площина якого проходить через деяку фіксовану пряму $l_{0}$ .
  Ясно, що множина $W$ гомеоморфна $\mathbb{R}P^{2}\times \mathbb{R}P^{1}$ .
  Якщо ми введемо в просторі $\mathbb{R}^{3}$ координати так,
   що  $$l_{0}=\{x=(x_{1},x_{2},x_{3})\in \mathbb{R}^{3}
| x_{1}=x_{2}=0\}$$  то кожна пряма   з множини $W$ буде
задаватися рівнянням
$$l=\{x=(x_{1},x_{2},x_{3})\in \mathbb{R}^{3} |
ax_{1}+bx_{2}=0,
\frac{bx_{1}-ax_{2}}{\sqrt{a^{2}+b^{2}}}c+dx_{3}+e=0\},$$
$(a:b)\in \mathbb{R}P^{1}, (c:d:e) \in \mathbb{R}P^{2} $,   (якщо
$\{ e_{1} = (1,0,0),e_{2} = (0,1,0) \}$ - ортонормований базис на
площині ${x_{3}=0}$, тоді $$\{
\frac{be_{1}-ae_{2}}{\sqrt{a^{2}+b^{2}}}, e_{3} \}$$ -
ортонормований базис на площині $ax_{1}+bx_{2}=0$ і однією з
координат на ній можна вважати число $\displaystyle
\frac{bx_{1}-ax_{2}}{\sqrt{a^{2}+b^{2}}}$).
 Це рівняння не є лінійним відносно координат спряженого простору,
  а точці $x=(0,0,x_{3})$  відповідає двовимірна площина.
   Але, якщо ми в множині $W$ замінимо
  прямі, паралельні площині $S$, площинами паралельними площині $S$,
   то рівняння $ax_{1}+bx_{2}=0$  ми можемо записати у вигляді
  $$c\frac{ax_{1}}{\sqrt{a^{2}+b^{2}}}+c\frac{bx_{2}}{\sqrt{a^{2}+b^{2}}}=0,$$
  яке перетворюється в тотожність при $c=0$. Замінивши
  $\displaystyle\frac{bc}{\sqrt{a^{2}+b^{2}}}$
   на $a'$, $\displaystyle\frac{-ac}{\sqrt{a^{2}+b^{2}}}$  на $b'$,
    $d$ на $c'$, $e$ на $d'$  одержимо
    $$\displaystyle l=\{x=(x_{1},x_{2},x_{3})\in \mathbb{R}^{3}
| a'x_{1}+b'x_{2}+c'x_{3}+d'=0, -b'x_{1}+a'x_{2}=0\},$$
$$(a')^{2}+(b')^{2}+(c')^{2}+(d')^{2}\neq 0.$$ Отже множина $W$
ізоморфна $\mathbb{R}P^{3}$ й складається із прямих
   паралельних даній прямій або тих, що перетинають її під гострим кутом, а
   також з перпендикулярних їй площин.
   Аналогічно попередньому прикладу, $\mathbb{R}^{3}$ можна
   розширити до $\mathbb{R}P^{3}$
   і тоді відношення спряження буде симетричним і лінійним
    відносно двох ізоморфних просторів.

    \vskip 5mm

\exam Розглянемо більш загальний приклад. Візьмемо лінійний простір
$\mathbb{R}P^{n}$ і кососиметричний оператор $A$ на ньому. Відомо,
що $(Ax,x)=0$, де дужки означають скалярний добуток. Нехай $W(A)$
складається з лінійних підпросторів, ортогональних векторам $Ax$ й
$x$ для деякого вектора $x$. Елементами $W(A)$ будуть площини
корозмірності 1 й 2. Якщо $A=0$, то множина $W(A)$ буде множиною
всіх гіперплощин. В останньому випадку кожна площина, що не
проходить через початок координат, задається точкою $y \in
\mathbb{R}P^{n}$ за допомогою спряження $(x,y)=0$. Помітимо, що
множина площин, що проходять через дану точку, задається точками,
які лежать на деякій площині. Цю відповідність можна узагальнити й
для довільного оператора $А$ за допомогою рівностей $(x,y)=0$,
$(x,Ay)=0$. Тоді множина, спряжена точці, задається рівностями
$(y,x)=0$, $(y,Ax)=0$  і належить $W(A)$ (оскільки
$(x,Ay)=-(Ax,y)=0$ для кососиметричного оператора). У зв'язку з тим,
що для непарних розмірностей всі лінійні оператори мають інваріантні
прямі, то множина $W(A)$ для непарних розмірностей обов'язково
містить гіперплощини (корозмірності 1) поряд із площинами
корозмірності~2.

 \vskip 0.5cm
\setcounter{theorem}{0} \setcounter{lemma}{0}
\addtocounter{section}{1} \setcounter{punkt}{0}
\section*{\thesection. \
Основна теорема} \addcontentsline{toc}{section}{\thesection. \
Основна теорема}

\theor \label{th} {\em Нехай $D$ відкрита обмежена зв'язна множина
в $\mathbb{R}^{n}$ із гладкою межею класу $\mathbb{C}^{1}$, опукла
відносно сім'ї $(2n-2)$-вимірних площин, паралельних площині
$\{x\in \mathbb{R}^{n}\mid x_{n}=0\}$. Тоді її межа $\partial D$ є
когомологічною сферою $S^{n-1}$.}

\textbf{Доведення.} Розглянемо проекцію $\lambda$ множини $Cl \ D$
на пряму $Ox_{n}$. Її звуження на межу $\lambda:
\partial D \rightarrow \mathbb{R}$ -
гладке відображення гладких многовидів. Відповідно до теореми
Сарда [4] образ множини його сингулярних точок $\lambda (S),
S\subset \partial D$ має міру нуль на прямій. Отже, для майже всіх
точок $y \in \lambda (\partial D)$ множина $\lambda^{-1} (y)
\bigcap \partial D$ є сукупністю гладких поверхонь, межею лінійно
опуклої відкритої множини в площині $\{x \in \mathbb{R} |
x_{n}=y\}=\Lambda_{y}$.

 Для регулярних значень $y$ відображення $\lambda$
 множини $\lambda^{-1}_{y} \subset Cl \ D$ і
 $Cl (\lambda^{-1}_{y} \bigcap D)$
 співпадають. А отже $\Lambda_{y} \bigcap \partial D$ - межа опуклої
  області $\lambda^{-1}(y)=\Lambda_{y} \bigcap Cl \ D$ у гіперплощині
$\Lambda_{y}$. Нехай
  $x_{0} \in \lambda^{-1}(\lambda(S))\bigcap D$  - довільна точка
нерегулярного перетину області $D$. Виберемо окіл $U(x_{0})$ цієї
точки, що лежить в $D$. Нехай $x_{1}$ інша точка множини
$\lambda^{-1}(\lambda(S))\bigcap Cl \ D$. У силу щільності в
$\lambda (Cl \ D)$ регулярних перетинів, існує послідовність
$A=\{x_{1n}\}$ точок, що належать регулярним перетинам, яка
сходиться до точки $x_{1}$. Виберемо підпослідовність
послідовності $A$, що має властивість $\lambda(x_{1n})\subset
\lambda(U(x_{0}))$. Не порушуючи загальності нехай цією
властивістю володіє сама послідовність $A$. Виберемо іншу
послідовність $B=\{x_{0n}\}$ збіжну до точки $x_{0}$, і таку, що
має додаткову властивість $x_{0n} \in U(x_{0}) \bigcap
\lambda^{-1}(\lambda(x_{1n})) $. Пари точок $x_{1n}$ й $x_{0n}$
належать опуклій множині $D \bigcap
\lambda^{-1}(\lambda(x_{1n}))$, і тому відрізок, що їх з'єднує,
$[x_{1n},x_{0n}]$ теж належить $D \bigcap
\lambda^{-1}(\lambda(x_{1n}))$. Границя послідовності відрізків
$[x_{1n},x_{0n}]$ природно містить у собі відрізок
$[x_{1},x_{0}]$. Тому для довільної точки $y \in \lambda(Cl \ D)$
множина $\lambda^{-1}(y)\bigcap Cl \ D$ буде зірковою відносно
довільної точки з $\lambda^{-1}(y)\bigcap D$, а множина
$\lambda^{-1}(y)\bigcap D$ опуклою. Аналогічно доводиться, що
множина $\lambda^{-1}(y)$ зіркова для межових точок множини
$\lambda(Cl \ D)$.

Отже, перетин будь-якого променя паралельного площині $\{x \in
\mathbb{R}^{n}| x_{n}=0 \}$, що виходить із області $D$, із межею
$\partial D$ є або точкою  або замкнутим відрізком (інакше кажучи,
він є ациклічним).

Доведемо, що існує крива $l \subset Cl \ D$, що взаємно однозначно
проектується на $\lambda(Cl \ D)$. Нехай точки $a,b \in \partial
D$ проектуються в кінці відрізка $\lambda(\bar{D})$. Внаслідок
зв'язності існує спрямлювана крива $l_{1}$, яка їх з'єднує і
повністю, за винятком своїх кінців, лежить в $D$. В околах точок
$a$ й $b$ криву $l_{1}$ задамо відрізками $[a,a_1]$, $[b,b_1]$
внутрішніх нормалей до гладкої межі так, щоб вони проектувалися
взаємно однозначно на свій образ. Будемо йти уздовж кривої до
точки $c$, для якої множина $\lambda^{-1}(\lambda(c))\bigcap
l_{1}$ складається з більш ніж однієї точки.

\begin{figure}[h!]
  \centering
  \includegraphics[width=200pt]{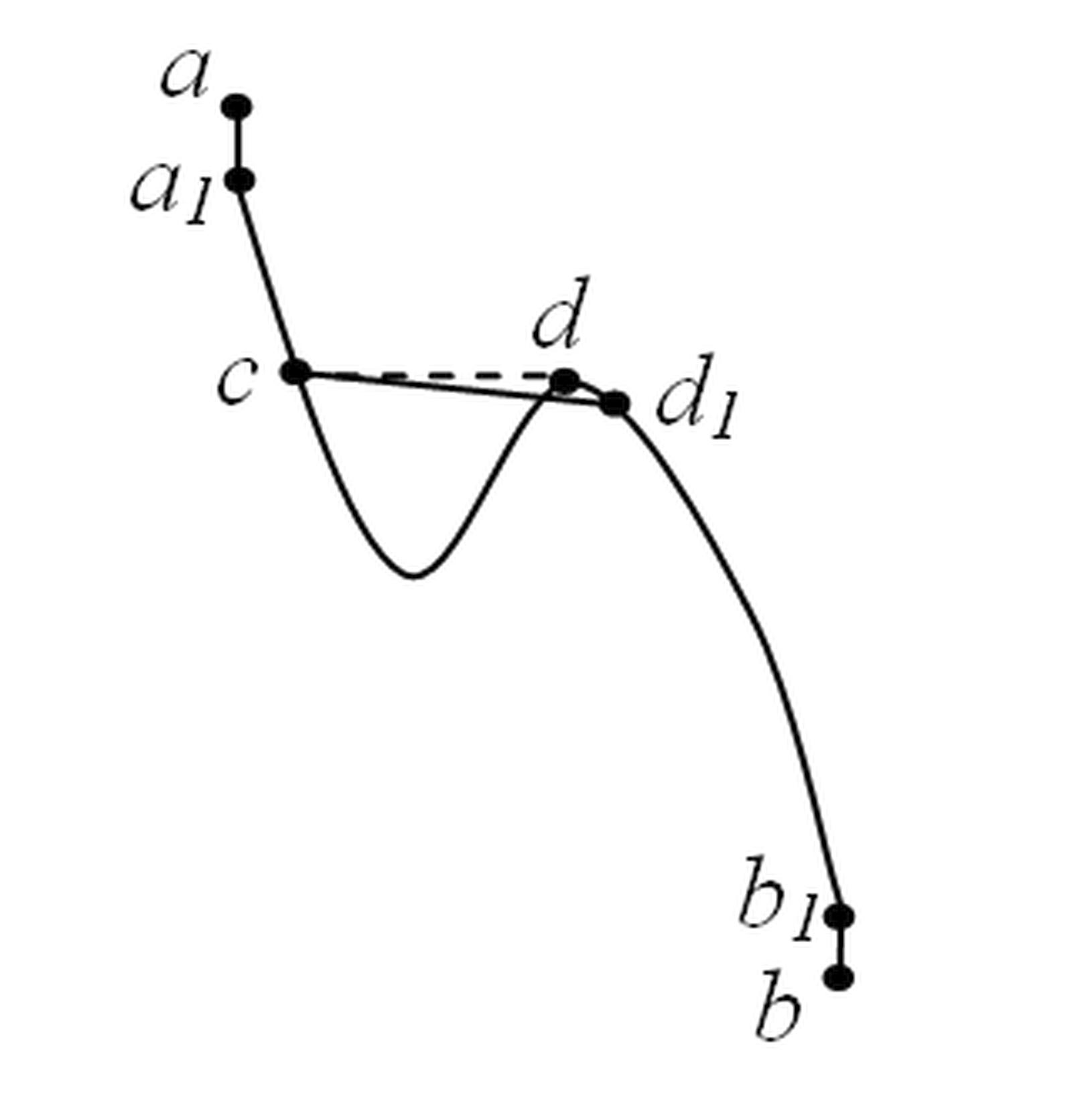}
\end{figure}

 Якщо ця точка не належить регулярному перетину, то
замість такої точки $с$ візьмемо точку $c_{1}$, що слідує за $c$
на кривій $l_{1}$ і таку, що $\lambda^{-1}(\lambda(c_{1}))\bigcap
l_{1}$ також складається з більш ніж однієї точки. Вибір такої
точки можливий у силу відкритості $D$. У цій множині візьмемо
точку $d$, яка на кривій $l_{1}$ слідує за всіма іншими точками
$\lambda^{-1}(\lambda(c))\bigcap l_{1}$, це можливо в силу
компактності кривої. Нехай $\Lambda_{\alpha} = \{x \in
\mathbb{R}^{n} | x_{n}=\alpha\}$. Відрізок $[c,d]$ повністю лежить
в $D$ у силу опуклості $D\bigcap \Lambda_{\lambda(c)}$ і
проектується в одну точку відрізка $\lambda(Cl \ D)$. Оскільки
область $D$ відкрита, існує точка $d_{1}$, що слідує за $d$ на
кривій, така, що відрізок $[c,d_{1}]$ лежить в $D$, а дуга кривої
$\smile ac$ в об'єднанні з відрізком $[c,d_{1}]$ проектується на
пряму $Ox_n$ взаємно однозначно. Оскільки крива $l_{1}$ є
впорядкованою множиною, то, пройшовши уздовж кривої і замінивши
послідовно необхідні дуги на відрізки, ми одержимо криву $l
\subset Cl \ D$, що взаємно однозначно проектується в пряму
$Ox_n$.

\begin{figure}[h!]
  \centering
  \includegraphics[width=300pt]{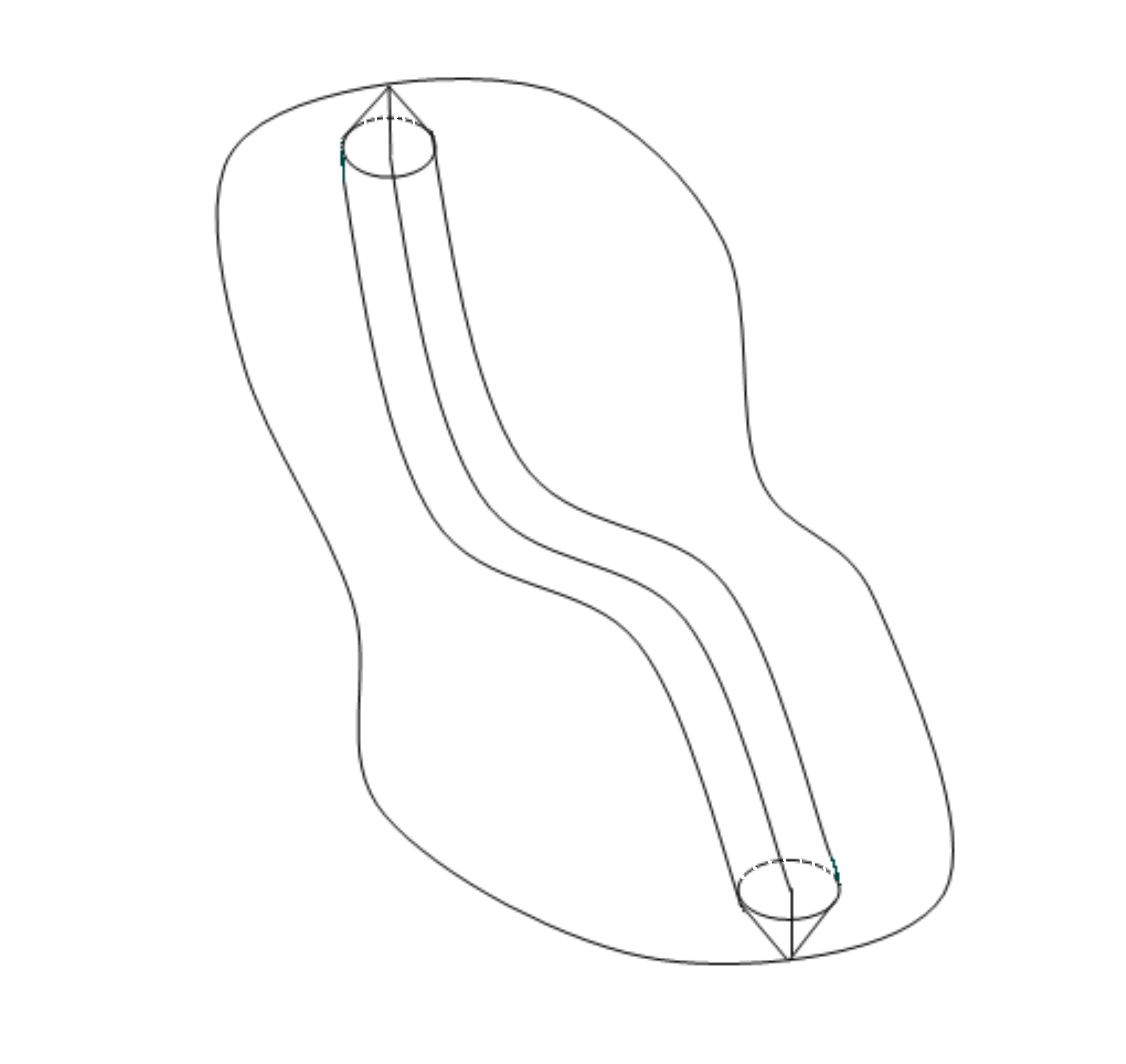}
\end{figure}

Нехай $\lambda(D)=(\alpha,\beta)$, тоді для як завгодно малого
$\varepsilon >0$ існує $\delta >0$ таке, що для будь-якого $\gamma
\in [\alpha+\varepsilon,\beta-\varepsilon]$ існує сфера в площині
$\Lambda_{\gamma}$ радіуса $\delta$ із центром у точці $l\bigcap
\Lambda_{\gamma}$, що цілком лежить в $D$. Якщо необхідно,
зменшимо $\varepsilon$ й $\delta$, так щоб конуси з вершинами в
точках $a$ й $b$ і основами $S^{n-2}(\alpha+\varepsilon)$ й
$S^{n-2}(\beta-\varepsilon)$ відповідно належали $D$. Об'єднання
$H$ всіх таких сфер і конусів гомеоморфне сфері $S^{n-1}$.

Розглянемо відображення $\sigma: \partial D \rightarrow H$, таке
що точка $x \in \partial D$ відображається в точку перетину
променя з початком в $\Lambda_{\lambda(x)}\bigcap l$, що проходить
через $x$, і сфери $S^{n-2}(\lambda(x))$.

 Таке відображення є однозначним й ациклічним (прообразом точки є або
точка або замкнутий відрізок). А отже, внаслідок теореми
Вієторіса-Бегля [\ref{spen}] воно індукує
ізоморфізм груп когомологій множин $\partial D$ й $S^{n-1}$.
Теорема доведена.$\blacksquare$

\textbf{Наслідок 1.} {\em Нехай $D$ відкрита обмежена зв'язна
множина в $\mathbb{R}^{n}$ із гладкою межею класу $\mathbb{C}^{1}$,
опукла відносно сім'ї $(2n-2)$-вимірних площин, паралельних площині
$\{x\in \mathbb{R}^{n}\mid x_{n}=0\}$. Тоді її межа $\partial D$ є
топологічною сферою $S^{n-1}$ при $n > 4$.}

\textbf{Доведення.} Розглянуте в доведенні теореми відображення
$\sigma: \partial D \rightarrow H$  є однозначним й точковим
(прообразом точки є або точка або замкнутий відрізок). А отже,
внаслідок теореми Компанійця [\ref{komp}] воно
індукує ізоморфізм груп гомотопій множин $\partial D$ й $S^{n-1}$.

Згідно Теореми Смейла про $h$-кобордизм [\ref{Smale}] та теореми
Фрідмана [\ref{Freedman}] при $n > 4$ гомотопічна сфера є
топологічною сферою.     $\blacksquare$

Нехай тепер $W_1$ -- множина площин в просторі $\mathbb{R}^n$
 таких, що афінна оболонка, натягнута на довільну
з них і на фіксовану площину $\Lambda_0$ корозмірності $2$, не
співпадає з $\mathbb{R}^n$. Іншими словами, $W_1$ -- множина площин,
які належать гіперплощинам пучка гіперплощин $\Delta_\alpha$, що
містять площину $\Lambda_0$.

 \textbf{Наслідок 2.} {\em Нехай $D$ відкрита обмежена зв'язна множина
в $\mathbb{R}^{n}$ із гладкою межею класу $\mathbb{C}^{1}$, опукла
відносно сім'ї $W_1$. Якщо $D$ не має спільних точок з деякою
площиною пучка $\Delta_\alpha$ (позначимо її $\Delta_0$), тоді її
межа $\partial D$ є когомологічною сферою $S^{n-1}$, а при $n
> 4$ межа $\partial D$ є також гомеоморфною сфері $S^{n-1}$.}

\textbf{Доведення.} Розглянемо довільне проективне відображення
$\pi$ простору $\mathbb{R}^n$ самого в себе, яке переводить площину
$\Delta_0$ в нескінченність. Це відображення переводить множину
$\mathbb{R}^n \backslash \Delta_0$ гомеоморфно саму в себе, а
гіперплощини пучка $\Delta_\alpha$ в паралельні гіперплощини. Отже,
множина $\pi(D)$ є пошарово лінійно опуклою, і внаслідок попередньої
\textbf{теореми 2.3.1} її межа $\partial(\pi(D))$ є когомологічною
сферою. Оскільки, гомеоморфні множини мають ізоморфні групи
когомологій, то і $\partial D$ є когомологічною сферою $S^{n-1}$.
Твердження про гомеоморфність $\partial D$ сфері $S^{n-1}$ при $n
> 4$ доводиться аналогічно попередньому наслідку.
$\blacksquare$

Наступні приклади показують істотність гладкості межі й
обмеженості області в теоремі.

\exam Нехай $$D=\{x=(x_{1},x_{2},x_{3})\in \mathbb{R}^{3} |
1<|x_{1}|+|x_{3}|<3, 0<x_{2}<1\}.$$ Легко переконатися, що область
$D$ опукла відносно сімейства прямих, паралельних площині
$\{x_{3}=0\}$, оскільки через кожну доповнення проходить пряма
паралельна або осі $0x_{1}$ або осі $0x_{2}$ і така, що не
перетинає області $D$. Перетини цієї області площинами
$\{x_{3}=d\}$ при $1 < |d| < 3$ будуть прямокутниками
$$\{x=(x_{1},x_{2},x_{3})\in \mathbb{R}^{3} | |x_{1}|<3-d,
0<x_{2}<1\},$$ а при $-1\leq d\leq 1$ будуть складатися із двох
прямокутників $$\{x=(x_{1},x_{2},x_{3})\in \mathbb{R}^{3} | (d-3 <
x_{1} < d-1, 0<x_{2}<1) \bigvee$$ $$\bigvee (1-d < x_{1} < 3-d,
0<x_{2}<1) \}.$$

 Границя області $\partial D$ гомеоморфна тору й
має нетривіальну одновимірну групу когомологій.

\begin{figure}[h!]
  \centering
  \includegraphics[width=300pt]{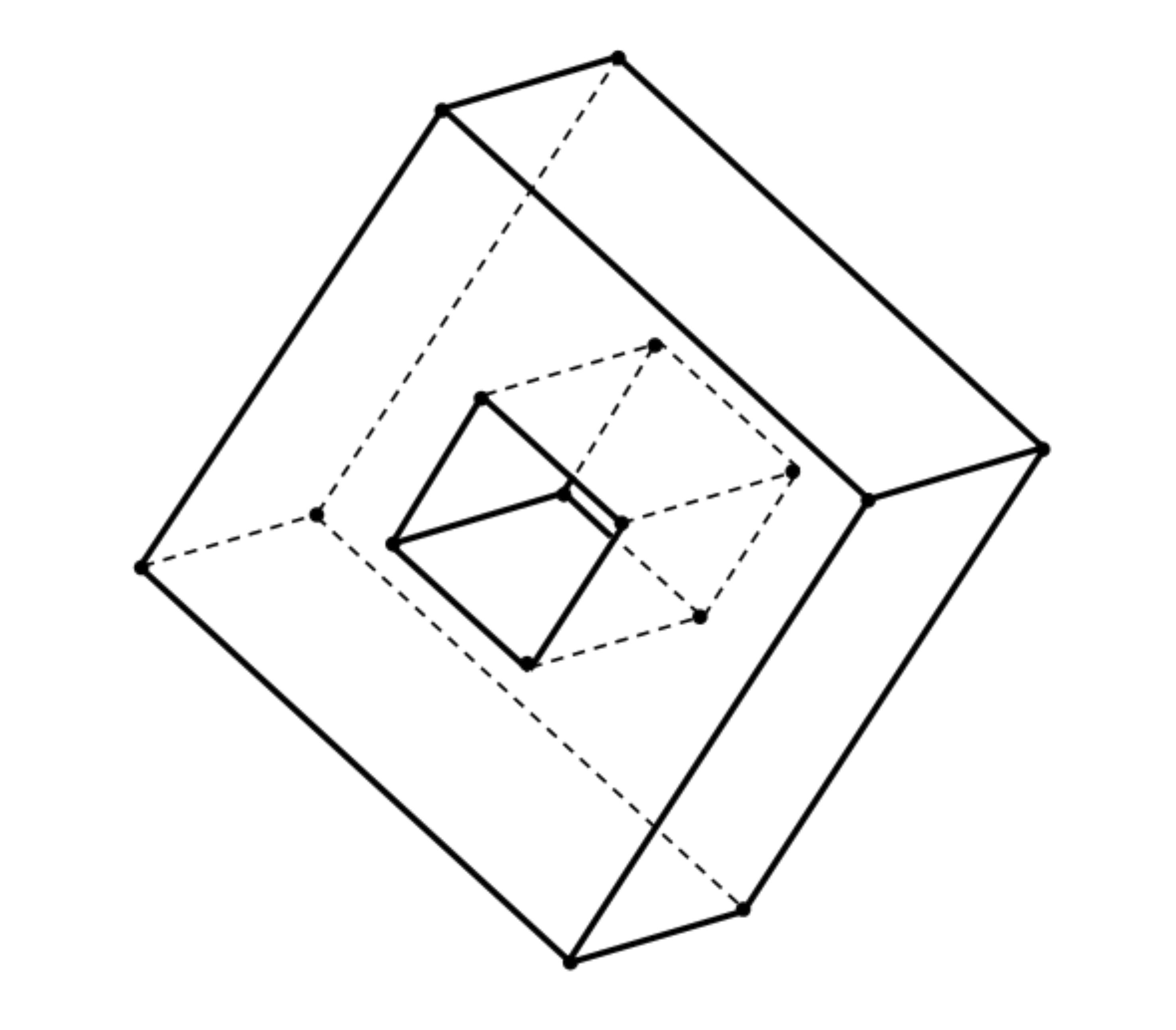}
\end{figure}

\exam Нехай $$D=\{x=(x_{1},x_{2},x_{3})\in \mathbb{R}^{3} |
f(x_{1},x_{2},x_{3})=x_{1}^{2}-x_{2}^{2}(1-x_{3}^2)+x_{3}^{2}-1>0\}.$$
Кожен перетин області $D$ площиною $$\Sigma_{d} = \{ x_{3}=a\},
-1<a<1$$ являє собою $$D \bigcap \Sigma_{d} = \{x=(x_{1},
x_{2})\in \mathbb{R}^{2} | x_{1}^{2}/(1-d^2)-x_{2}^{2}>1\}$$
внутрішність гіперболи (дві компоненти, що не перетинають пряму
${x_1=0}$) і є лінійно опуклим. Область $D$ має гладку межу,
оскільки градіент
$$grad \ f = (2x_{1},2x_{2}(1-x_{3}^2), 2x_{3}(1-x_{2}^2))$$
відмінний від 0 на $\partial D$, і опукла відносно сім'ї прямих,
паралельних площині $\{x_{3}=0\}$. Межа області $\partial D$
гомеоморфна поверхні нескінченного циліндра над колом й має
нетривіальну одновимірну групу когомологій.

\begin{figure}[h!]
  \centering
  \includegraphics[width=300pt]{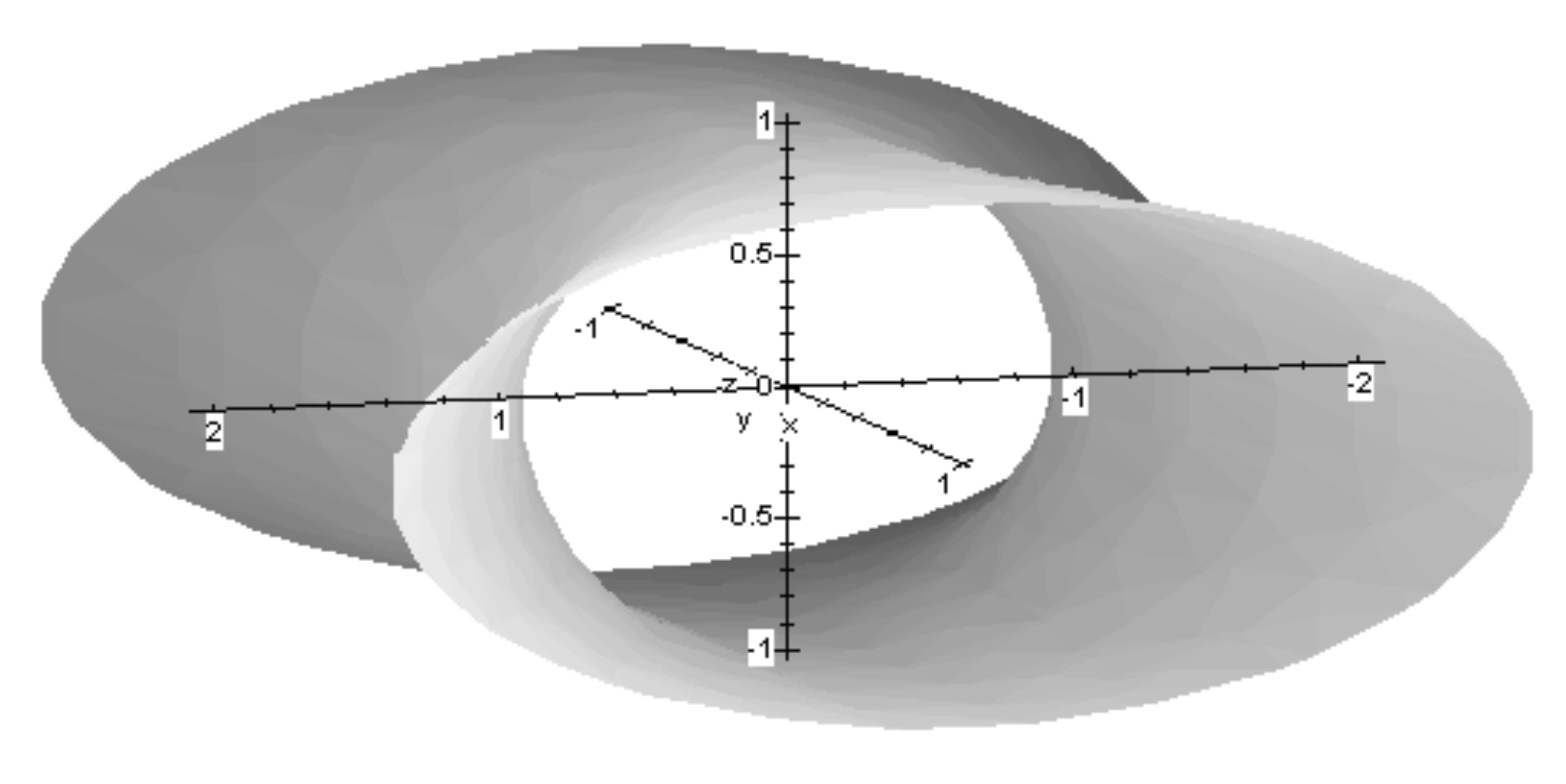}
\end{figure}


\end{document}